\documentclass[
leqno,          
12pt,           
a4paper         
]{amsart}

\usepackage[colorlinks,citecolor=blue,urlcolor=blue]{hyperref}          
\usepackage{courier}                                                    
\usepackage{amssymb, amsmath, amsfonts, amsthm}                         
\usepackage{mathtools, cite, enumerate, rotating, framed, lipsum}       
\usepackage{fancybox, bbm}                                              
\usepackage[dvipsnames]{xcolor}                                         
\usepackage[polish, english]{babel}                                     
\usepackage[T1]{fontenc}                                                
\usepackage[cp1250]{inputenc}                                           
\usepackage[normalem]{ulem}
\usepackage{yfonts}
\usepackage{mathrsfs}

\usepackage[a4paper,top=3cm,bottom=2cm,left=3cm,right=3cm,marginparwidth=1.75cm]{geometry}

\usepackage{amsmath,amstext,amssymb,amsopn,amsthm,enumitem}
\usepackage{graphicx}
\usepackage[colorinlistoftodos]{todonotes}
\usepackage{yfonts}
\usepackage{mathrsfs}
\usepackage{xstring}

\DeclareMathAlphabet{\mathpzc}{OT1}{pzc}{m}{it}
 \numberwithin{equation}{section}                        

\newcommand{\thmcount}{equation}                 
\newcounter{specialcounter}

\newtheorem{Thm}[\thmcount]{Theorem}
\newtheorem{Sthm}[specialcounter]{Theorem}
\newtheorem{Cor}[\thmcount]{Corollary}
\newtheorem{Lem}[\thmcount]{Lemma}
\newtheorem{Prop}[\thmcount]{Proposition}
\newtheorem{Rem}[\thmcount]{Remark}
\newtheorem{Defn}[\thmcount]{Definition}
\newtheorem{Ex}[\thmcount]{Example}
\newtheorem{Asu}[\thmcount]{Assumption}
\newtheorem{Sol}[\thmcount]{Solution}
\newtheorem*{Thmx}{Theorem}
\newtheorem*{Corx}{Corollary}
\newtheorem*{Lemx}{Lemma}
\newtheorem*{Propx}{Proposition}
\newtheorem*{Remx}{Remark}
\newtheorem*{Defnx}{Definition}
\newtheorem*{Exx}{Example}
\newtheorem*{Asux}{Assumption}
\newtheorem*{Solx}{Solution}

\newcommand \eq[1]{\begin{equation} #1 \end{equation}}
\newcommand \eqx[1]{\begin{equation*}  #1 \end{equation*}}
\newcommand \al[1]{\begin{align} #1 \end{align}}
\newcommand \alx[1]{\begin{align*}  #1 \end{align*}}
\renewcommand \sp[1]{\begin{equation} \begin{split} #1 \end{split} \end{equation}}
\newcommand \spx[1]{\begin{equation*} \begin{split} #1 \end{split} \end{equation*}}
\newcommand \en[1]{\begin{enumerate}  #1 \end{enumerate}}
\newcommand \ite[1]{\begin{itemize}  #1 \end{itemize}}

\newcommand{\sthm}[2]{\begin{Sthm} \label{#1} #2 \end{Sthm}}

\newcommand{\lem}[2]{\begin{Lem} \label{#1} #2 \end{Lem}}

\newcommand{\cor}[2]{\begin{Cor} \label{#1} #2 \end{Cor}}

\newcommand{\prop}[2]{\begin{Prop} \label{#1} #2 \end{Prop}}

\newcommand{\rem}[2]{\begin{Rem} \label{#1} #2 \end{Rem}}

\newcommand{\defn}[2]{\begin{Defn} \label{#1} #2 \end{Defn}}

\newcommand{\pr}[1]{\begin{proof} #1 \end{proof}}

\newcounter{comcount}

\renewcommand{\hline}{\vbox{\hrule width\textwidth height 1pt}\smallskip}

\renewcommand{\a}{\alpha}
\newcommand{\be}{\beta}         \newcommand{\e}{\varepsilon}
                 \newcommand{\de}{\delta}
        \newcommand{\la}{\lambda}

\newcommand{\abs}[1]{\left| #1 \right|}
\newcommand{\set}[1]{\left\{ #1 \right\}}
\newcommand{\norm}[1]{\left\| #1 \right\|}

\newcommand{\eee}[1]{\left( #1 \right)}

\renewcommand{\al}{\alpha}

\newcommand{\wt}{\widetilde}
\newcommand{\LLag}{L^{[\a]}_{L}}

\newcommand{\CC}{\mathbb{C}}

\newcommand{\NN}{\mathbb{N}}

  \newcommand{\qq}{\mathcal{Q}}
\newcommand{\RR}{\mathbb{R}}  
  
 \newcommand{\Tt}{\mathbf{T}}

\newcommand{\ZZ}{\mathbb{Z}}  

\newcommand{\supp}{\mathrm{supp}}
\newcommand{\8}{\infty}

\newcommand{\Rd}{{\RR^d}}

\renewcommand{\rm}[1]{\mathrm{#1}}

\renewcommand{\Tt}{T_t}
\newcommand{\tti}[1]{T^{[#1]}}

\newcommand{\s}[2]{%
    \IfEqCase{#2}{%
        {1}{#1^*}%
        {2}{#1^{**}}%
        {3}{#1^{***}}
    }
}
\newcommand{\Ht}{H_t}

\title{Local atomic decompositions for multidimensional Hardy spaces }

\author[ Edyta Kania-Strojec ]{ Edyta Kania-Strojec }
 \address{
 Edyta Kania-Strojec \newline
 \indent Instytut Matematyczny, Uniwersytet Wroc\l awski \newline
 \indent pl. Grunwaldzki 2/4, 50-384 Wroc\l aw, Poland }
 \email{edyta.kania-strojec@math.uni.wroc.pl }

 \author[ Pawe\l \ Plewa ]{ Pawe\l \ Plewa }
 \address{
 Pawe\l \ Plewa \newline
 \indent Faculty of Pure and Applied Mathematics, Wroc\l aw University of Science and Technology \newline
 \indent Wyb. Wyspia\'nskiego 27, 50-370 Wroc\l aw, Poland }
 \email{pawel.plewa@pwr.edu.pl}

 \author[ Marcin Preisner ]{ Marcin Preisner }
 \address{
 Marcin Preisner \newline
 \indent Instytut Matematyczny, Uniwersytet Wroc\l awski \newline
 \indent pl. Grunwaldzki 2/4, 50-384 Wroc\l aw, Poland }
 \email{marcin.preisner@uwr.edu.pl }

\thanks{ The first and the third authors are supported by the grant No. 2017/25/B/ST1/00599 from National Science Centre (Narodowe Centrum Nauki), Poland. The second author is supported by the National Science Centre of Poland, NCN grant No. 2018/29/N/ST1/02424. }

\begin{document}
\maketitle

\begin{abstract}
We consider a nonnegative self-adjoint operator $L$ on $L^2(X)$, where $X\subseteq \mathbb{R}^d$. Under certain assumptions, we prove atomic characterizations of the Hardy space
$$H^1(L) = \set{f\in L^1(X) \ : \ \norm{\sup_{t>0} \abs{\exp(-tL)f }}_{L^1(X)}<\8}.$$
We state simple conditions, such that $H^1(L)$ is characterized by atoms being either the classical atoms on $X\subseteq \Rd$ or local atoms of the form $|Q|^{-1}\chi_Q$, where $Q\subseteq X$ is a cube (or cuboid).

One of our main motivation is to study multidimensional operators related to orthogonal expansions. We prove that if two operators $L_1, L_2$ satisfy the assumptions of our theorem, then the sum $L_1 + L_2$ also does. As a consequence, we give atomic characterizations for multidimensional Bessel, Laguerre, and Schr\"odinger operators.

As a by-product, under the same assumptions, we characterize $H^1(L)$ also by the maximal operator related to the subordinate semigroup $\exp(-tL^\nu)$, where $\nu\in(0,1)$.
\end{abstract}

\section{Background and main results}\label{sec1}

\subsection{Introduction}\label{ssec11}
Let us first recall that the classical Hardy space $H^1(\Rd)$ can be defined by the maximal operator, i.e.
$$f\in H^1(\Rd) \qquad \iff \qquad \sup_{t>0} \abs{H_t f} \in L^1(\Rd). $$
Here and thereafter $H_t = \exp(t\Delta)$ is the heat semigroup on $\Rd$ given by $H_t f(x) = \int_{\Rd} H_t(x,y) f(y)  \, dy$,
\eq{\label{heat_kernel}
H_t(x,y) = (4\pi t)^{-d/2} \exp\eee{-\frac{|x-y|^2}{4t}}.
}

Among many  equivalent characterizations of $H^1(\Rd)$ one of the most useful is the characterization by atomic decompositions proved by Coifman \cite{Coifman_Studia} in the one-dimensional case and by Latter \cite{Latter_Studia} in the general case $d\in \NN$. It says that $f\in H^1(\Rd)$ if and only if $f(x) = \sum_{k=1}^\8 \la_k a_k(x)$, where $\la_k \in \CC$ are such that $\sum_{k=1}^\8 \abs{\la_k} <\8$ and $a_k$ are {\it atoms}. By definition, a function $a$ is {\it an atom} if there exists a ball $B \subseteq \Rd$ such that:
\eqx{\label{atoms_class}
\supp \, a \subseteq B, \qquad \norm{a}_\8 \leq |B|^{-1}, \qquad \int_B a(x) \, dx = 0,
}
i.e. $a$ satisfies well-known localization, size, and cancellation conditions.

Later, Goldberg in \cite{Goldberg_Duke} noticed that if we restrict the supremum in the maximal operator above to the range $t\in(0, \tau^2)$, with $\tau>0$ fixed, then still the atomic characterization holds, but with additional atoms of the form $a(x) = |B|^{-1} \chi_B(x)$, where $\chi$ is the characteristic function and $B$ is a ball of radius~$\tau$ (see Section \ref{sec_local} for details).

Then, many atomic characterizations were proved for various operators including operators with Gaussian (or Davies-Gaffney) estimates, operators on spaces of homogeneous type, operators related to orthogonal expansions, Schr\"odinger operators, and others. The reader is referred to \cite{Hofmann_Memoirs, DZ_Studia_DK, Dziubanski2017,Dziubanski_Houston,DPRS,BDT_d'Analyse, Song_Yan, Auscher_unpublished,Stein} and references therein.

In this paper we deal with atomic characterizations of the Hardy space $H^1$ for operators, such that $H^1$ admits atoms of local type, i.e. atoms of the form $|B|^{-1} \chi_B$. We shall consider operators defined on $L^2(X)$, where $ X \subseteq \Rd$ with the Lebesgue measure. Our main focus will be on sums of the form $L=L_1+...+L_d$, where each $L_i$ acts only on the variable $x_i$, where $x=(x_1,...,x_d)$. For such $L$ we look for atomic decompositions. As an application, we can take operators related to some multidimensional orthogonal expansions. Additionally we prove characterizations of $H^1$ by subordinate semigroups.

\subsection{Notation}\label{ssec12}
Let $X = (a_1, b_1) \times ... \times (a_d,b_d)$ be a subset of $\Rd$. We allow $a_j=-\8$ and $b_j = \8$ so that we consider products of lines, half-lines, and finite intervals. We equip $X$ with the Euclidean metric and the Lebesgue measure. In the product case it is more convenient to use cubes and cuboids instead of balls, so denote for $z = (z_1, ..., z_d)\in X$ and $r_1,...,r_d>0$ the closed cuboid
$$Q(z,r_1,...,r_d) = \set{x\in X \ : \ |x_i-z_i|\leq r_i \text{ for } i=1,...,d},$$
and the cube
$Q(z,r) = Q(z,r,...,r).$
We shall call such $z$ the center of a cube/cuboid. For a cuboid $Q$ by $d_Q$ we shall denote the diameter of $Q$.

\defn{def_covering}{
Let $\qq$ be a  set of cuboids in $X$. We call $\qq$ {\it an admissible covering} of $X$ if there exist $C_1, C_2 >0$ such that:
\en{
\item $X = \bigcup_{Q\in \qq} Q$,
\item if $Q_1, Q_2 \in \qq$ and $Q_1\neq Q_2$ , then $|Q_1\cap Q_2| = 0$,
\item if $Q = Q(z,r_1,...,r_d)\in \qq$, then $r_i\leq C_1 r_j$ for $i,j \in \set{1,...,d}$,
\item if $Q_1, Q_2 \in \qq$ and $Q_1 \cap Q_2 \neq \emptyset$, then $C_2^{-1} d_{Q_1} \leq d_{Q_2}\leq C_2 d_{Q_1}$.}}

Let us note that {\bf 3.} means that our cuboids are almost cubes. In fact, we shall often use only cubes.

By $Q^*$ we shall denote a slight enlargement of $Q$. More precisely, if $Q=(z,r_1,...,r_d)$, then $Q^*:= Q(z, \kappa r_1,...,\kappa r_d)$, where $\kappa>1$. Observe that if $\qq$ is admissible covering of $\RR^d$, then choosing $\kappa$ close enough to 1 the family $\set{Q^{***}}_{Q\in \qq}$ is a finite covering of $\Rd$, namely
\eq{\label{finite_covering}
\sum_{Q\in \qq} \chi_{Q^{***}}(x) \leq C, \qquad x\in\RR^d
}
and, for $Q_1,Q_2 \in \qq$,
\eq{\label{neighbours}
Q_1^{***}\cap Q_2^{***} \neq \emptyset \quad \iff \quad  Q_1 \cap Q_2 \neq \emptyset.
}
In this paper we always choose $\kappa$ such that \eqref{finite_covering} and \eqref{neighbours} are satisfied. Let us emphasize that $Q$ and $Q^{*}$ are always defined as a subset of $X$, not as a subset of $\Rd$.

Having two admissible coverings $\qq_1$ and $\qq_2$ on $\RR^{d_1}$ and $\RR^{d_2}$ we would like to produce an admissible covering on $\RR^{d_1+d_2}$. However, one simply observe that products $\set{Q_1 \times Q_2 \ : \ Q_1 \in \qq_1, Q_2 \in \qq_2}$, would not produce admissible covering (in general, {\bf 3.} would fail). Therefore, for the sake of this paper, let us state the following definition.

\defn{def_box_prod}{
Assume that $\qq_1$ and $\qq_2$ are admissible coverings of $X_1 \subseteq \RR^{d_1}$ and $X_2 \subseteq\RR^{d_2}$, respectively. We define an admissible covering of $X_1 \times X_2$ in the following way. First, consider the covering $\set{Q_1 \times Q_2 \ : \ Q_1 \in \qq_1, Q_2 \in \qq_2}$. Then we further split each $Q=Q_1 \times Q_2$. Without loss of generality let us assume that $d_{Q_1}>d_{Q_2}$. We split $Q_1$ into cuboids $Q_1^{[j]}$, $j=1,...,M$, such that all of them have diameters comparable to $d_{Q_2}$ and satisfy {\bf 3.} of Definition \ref{def_covering}. Then the cuboids $Q^{[j]} = Q_1^{[j]} \times Q_2$, $j=1,...,M$, satisfy:
\ite{
\item $Q = \bigcup_{j=1}^M Q^{[j]}$,
\item for $i,j \in \set{1,...,M}$, $i\neq j$, we have $|Q^{[i]} \cap Q^{[j]}| = 0$,
\item each $Q^{[j]}$ satisfies {\bf 3.} from Definition \ref{def_covering}.
}
Notice that $M\leq [d_{Q_1}/d_{Q_2}]^{d_1}$. We shall denote such covering by $\qq_1 \boxtimes \qq_2$.
}
One may check that the definition above leads to an admissible covering of $X_1\times X_2$.

Having an admissible covering $\qq$ of $X \subseteq \Rd$ we define a local atomic Hardy space $H^1_{at}(\qq)$ related to $\qq$ in the following way. We say that a function $a \, : \, X \to \CC$ is a $\qq-atom$ if:

\begin{itemize}
\item[(i)] either there is  $Q\in\qq$ and a cube $K\subset \s{Q}{1}$, such that:  $$\supp \, a \subseteq K, \ \ \norm{a}_\8 \leq |K|^{-1},  \ \ \int a(x)\, dx = 0;$$
\item[(ii)] or there exists $Q\in\qq$ such that $$\a(x) = |Q|^{-1}\chi_Q(x).$$
\end{itemize}

Having $\qq$-atoms we define {\it the local atomic Hardy space related to $\qq$, } $H^1_{at}(\qq)$, in a standard way. Namely, we say that a function $f$ is in $H^1_{at}(\qq)$ if $f(x) = \sum_k \la_k a_k(x)$ with $\sum_k |\la_k| <\8$ and $a_k$ being $\qq$-atoms. Moreover, the norm of $H^1_{at}(\qq)$ is given by
$$\norm{f}_{H^1_{at}(\qq)} = \inf \sum_k \abs{\la_k},$$
where the infimum is taken over all possible representations of $f(x) = \sum_k \la_k a_k(x)$ as above. One may simply check that $H^1_{at}(\qq)$ is a Banach space.

\label{ssec12}
In the whole paper by $L$ we shall denote a nonnegative self-adjoint operator  and by $\Tt = \exp(-tL)$ the heat semigroup generated by $L$. We shall always assume that there exists a nonnegative integral kernel $\Tt(x,y)$ such that $\Tt f(x) = \int_X \Tt(x,y) f(y) \, dy$.
Our initial definition of the Hardy space $H^1(L)$ shall be by means of the maximal operator associated with $\Tt$, namely
\eqx{ \label{max-hardy}
	H^1(L) = \set{ f\in L^1(X) \ : \ \norm{f}_{H^1(L)} := \norm{\sup_{t>0} \abs{ \Tt f}}_{L^1(X)}  < \8}.
   }
Moreover, we shall consider the subordinate semigroup $K_{t,\nu}=\exp(-tL^\nu)$, $\nu\in(0,1)$,  and its Hardy space, which is defined by
\eqx{ \label{max-hardy-sub}
	H^1(L^\nu) = \set{ f\in L^1(X) \ : \ \norm{f}_{H^1(L^\nu)} := \norm{\sup_{t>0} \abs{ K_{t,\nu} f}}_{L^1(X)}  < \8}.
   }

\subsection{Main results} \label{ssec13}

Let us assume that an admissible covering $\qq$ of $X$ is given. Recall that $H_t(x,y)$ is the classical semigroup on $\Rd$ given in \eqref{heat_kernel}, and denote by $P_{t,\nu}=\exp(-t(-\Delta)^\nu)$ the semigroup generated by $(-\Delta)^\nu$, $\nu\in(0,1)$, and given by $P_{t,\nu}f(x) = \int_{\Rd} P_{t,\nu}(x,y) f(y) \, dy$. The kernel $P_{t,\nu}(x,y)$ is a transition density of the symmetric $2\nu$-stable L\'evy process in $\Rd$. It is well-known that \eq{\label{poissone_kernel}
0\leq P_{t,\nu}(x,y) \leq C_{d,\nu} \ \frac{t}{\left(t^{1/\nu} + |x-y|^2\right)^{\frac{d}{2}+\nu}} ,\quad x,y \in \Rd, t>0, \ \nu\in(0,1),
}
see e.g. \cite[Subsec. 2.6]{MK_FCAA}, \cite{PG_WH_Studia}. Let us mention that in the particular case of $\nu = 1/2$, the semigroup $P_{t,1/2}$ is the well-known Poisson semigroup on $\RR^d$.

Assume that an operator $L$ is as in Subsection \ref{ssec12}. Let $\nu\in(0,1)$ and suppose that $\wt{T_t}(x,y)$ is either $H_t(x,y)$ or $P_{t^\nu, \nu}(x,y)$. Consider the following assumptions:
\begin{align}
\label{a0} \tag{$A_0'$}
&&&0\leq \Tt(x,y) \leq C \ \frac{t^\nu}{\left(t + |x-y|^2\right)^{\frac{d}{2}+\nu}}, &&x,y\in X,\ t>0,&&\\
 \label{a1} \tag{$A_1'$}
&&&\sup_{y \in \s{Q}{1}} \int_{(\s{Q}{2})^c} \sup_{t > 0} \Tt(x,y) dx \leq C,  &&Q\in \qq,
&&\\
\label{a2} \tag{$A_2'$}
&&&\sup_{y \in \s{Q}{1}} \int_{\s{Q}{2}} \sup_{t \leq  d_Q^2 } \abs{\Tt(x,y) - \wt{T_t}(x,y)}dx \leq C, &&Q\in \qq.
&&
\end{align}

\sthm{atom-hardy}{
Assume that for $L$, $T_t$, and an admissible covering $\qq$ the conditions \eqref{a0}--\,\eqref{a2} hold. Then $H^1(L) = H^1_{at}(\qq)$ and the corresponding norms are equivalent.}

The proof of Theorem \ref{atom-hardy} is standard and uses only local characterization of Hardy spaces as in \cite{Goldberg_Duke}. For the convenience of the reader we present the proof in Section \ref{sec3}.

Our first main goal is to describe atomic characterizations for sums of the form $L_1 + ... + L_N$, where each $L_j$ satisfies \eqref{a0}--\eqref{a2} on a proper subspace. This is very useful in many cases such as multidimensional orthogonal expansions. Instead of dealing with products of kernels of semigroups, we can consider only one-dimensional kernel, but we shall need to prove slightly stronger conditions. More precisely, we consider $X_1 \times ... \times X_N \subseteq \RR^{d_1} \times ... \times \RR^{d_N} = \Rd$. Assume that $L_i$ is an operator on $L^2(X_i)$, as in Subsection \ref{ssec12}. Slightly abusing the notation we keep the symbol $L_i$ for $I\otimes ...  \otimes L_i \otimes  ... \otimes I$ as the operator on $L^2(X)$ and denote
\eq{\label{Lsum}
Lf(x) = L_1f(x) + ... + L_Nf(x), \quad x = (x_1,...,x_N) \in X.
}

For $x_i,y_i \in X_i$, by $\tti{i}_t(x_i,y_i)$ we denote the kernel of $\tti{i}_t = \exp\eee{-tL_i}$. We shall assume that each $\tti{i}_t(x_i,y_i)$ , $i=1,...,N$, is nonnegative and has the upper Gaussian estimates, namely
 \eq{ \label{gauss} \tag{$A_0$}
 0\leq \tti{i}_t(x_i,y_i) \leq C_i t^{-d_i/2} \exp\left(-\frac{|x_i-y_i|^2}{c_it}\right), \quad x_i,y_i \in X_i, t>0.
} 	
Obviously, \eqref{gauss} implies \eqref{a0} for $T_t(x,y) = T_t^{[1]}(x_1,y_1)... T_t^{[N]}(x_N,y_N)$. Moreover, we shall assume that for each $i \in \set{1,...,N}$ there exist a proper covering $\qq_i$ of $\RR^{d_i}$ such that the following generalizations of \eqref{a1} and \eqref{a2} hold: there exists $\gamma \in (0,1/3)$ such that for every $\delta \in[0,\gamma)$ and every $i=1,..,N,$

\eq{ \label{a11} \tag{$A_1$}
\sup_{y \in Q^{*}} \int_{(Q^{**})^c} \sup_{t > 0 } t^\delta \tti{i}_t(x,y) dx \leq C d_Q^{2\delta}, \quad Q \in \qq_i,
}

\eq{ \label{a22} \tag{$A_2$}
\sup_{y \in Q^*} \int_{Q^{**}} \sup_{ t< d_Q^2} t^{-\delta}\abs{\tti{i}_t(x,y) - \Ht(x,y)}dx \leq C d_Q^{-2\delta} \quad Q \in \qq_i.
}
Here $H_t$ is the classical heat semigroup on $\RR^{d_i}$, depending on the context. Now, we are ready to state our first main theorem.

\sthm{mainthm}{
Assume that for $i=1,...,N$ kernels $\tti{i}_t(x_i, y_i)$ are related to $L_i$ and suppose that for $\tti{i}_t(x_i, y_i)$ together with admissible coverings $\qq_i$ the conditions \eqref{gauss}--\,\eqref{a22} hold. If $L=L_1+...+L_N$ is as in \eqref{Lsum}, then
\eqx{\label{hardy-equ}
H^1(L) = H^1_{at}(\qq_1 \boxtimes ... \boxtimes \qq_N)
}
and the corresponding norms are equivalent.
}

Our second main goal is to characterize $H^1(L)$ by the subordinate semigroup $K_{t,\nu} = \exp(-tL^\nu)$, for $0<\nu<1$. Obviously, one can try to apply Theorem \ref{atom-hardy}, but for many operators the subordinate kernel $K_{t,\nu}(x,y) $ is harder to analyze than $T_t(x,y)$ (e.g., in some cases a concrete formula with special functions exists for $T_t(x,y)$, but not for $K_{t,\nu}(x,y)$). However, it appears that under our assumptions \eqref{gauss}--\,\eqref{a22} we obtain the characterization by the subordinate semigroup essentially for free.

\sthm{mainthm2}{
Under the assumptions of Theorem \ref{mainthm}, for $\nu\in(0,1)$, we have that
\eqx{
H^1(L^\nu) = H^1_{at}(\qq_1 \boxtimes ... \boxtimes \qq_N).
}
Moreover, the corresponding norms are equivalent.
}

\subsection{Applications}\label{ssec13}

One of the goals of this paper is to verify the assumptions of Theorems \ref{mainthm} and \ref{mainthm2} for various well-known operators. In this subsection we provide the list of these operators.

\newcommand{\lbl}{L_B^{[\be]}}
\subsubsection{ {\bf Bessel operator.}}\label{sssec131}

For $\beta>0$ let $\lbl=-\frac{d^2}{dx} + \frac{\be^2-\be}{x^2}$ denote the one-dimensional Bessel operator on the positive half-line $X = (0,\8)$ equipped with the Lebesgue measure. The semigroup $T_{B,t}= \exp(-t \lbl)$ is given by $T_{B,t} f(x) = \int_X T_{B,t}(x,y)f(y)\, dy$, where
\eq{\label{bessel-kernel}
T_{B,t}(x,y) = \frac{(xy)^{1/2}}{2t} I_{\be-1/2}\left(\frac{xy}{2t} \right) \exp\left(-\frac{x^2+y^2}{4t} \right), \qquad x,y\in X, t>0.
}
Here, $I_\tau$ is the modified Bessel function of the first kind. The Hardy space $H^1(\lbl)$ for the one-dimensional Bessel operator was studied in \cite{BDT_d'Analyse}. In Section \ref{ssec41} we check that the assumptions \eqref{gauss}--\eqref{a22} are satisfied for $L_B$ with the admissible covering
$$\qq_B = \set{[2^n,2^{n+1}] \, : n\in \ZZ}$$
of $X=(0,\8)$. This gives a slightly simpler proof of the characterizations of $H^1(\lbl)$ by the maximal operators of the semigroups $\exp(-t\lbl)$ and, also, gives a characterization by $\exp(-t(\lbl)^\nu)$, $0<\nu<1$. We have the following corollary for the multidimensional Bessel operator.

\cor{coro_bess}{
Let $\be_1,...,\be_d>0$ and $L_B = L_B^{[\be_1]}+...+ L_B^{[\be_d]} $, be the multidimensional Bessel operator on $L^2((0,\infty)^d)$. Then, the Hardy spaces $H^1(L_B)$, $H^1(L_B^\nu)$, $\nu\in(0,1)$, and $H^1_{at}(\qq_B \boxtimes ... \boxtimes \qq_B)$ coincide. Moreover, the associated norms are comparable.
}


\begin{figure}[h]
\caption{The covering  $\qq_B \boxtimes \qq_B$}
\centering
\includegraphics[scale=0.45]{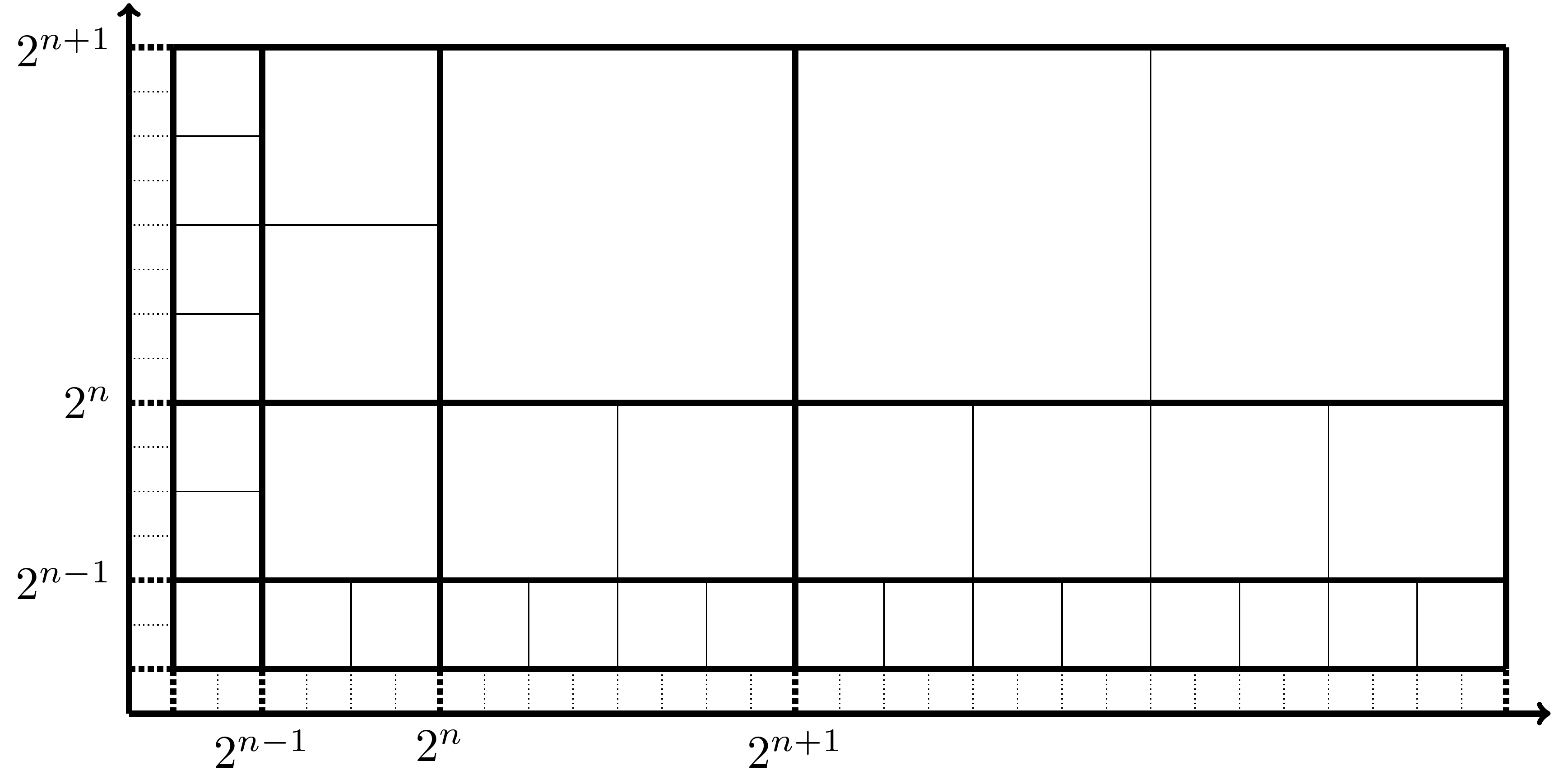}
\end{figure}

\subsubsection{ {\bf Laguerre operator.}}\label{sssec133}
Let $\a>-1/2$ and $L_{L}^{[\a]}=-\frac{d^2}{dx^2}+x^2+\frac{\al^2-1/4}{x^2}$ denote the Laguerre operator on $X=(0,\infty)$. The kernels associated with the heat semigroup $T_{L,t}=\exp\left(-tL_{L}^{[\a]}\right)$ are defined by
\eq{\label{lag1-kernel}
T_{L,t}(x,y) = \frac{(xy)^{1/2}}{\sinh 2t} I_{\a}\left(\frac{xy}{\sinh 2t} \right) \exp\left(-\frac{\cosh 2t}{2\,{\sinh{2t}}}(x^2+y^2) \right), \qquad x,y\in X, t>0.
}

\begin{figure}[h]
\caption{The covering  $\qq_{L} \boxtimes \qq_{L}$}
\centering
\includegraphics[scale=0.45]{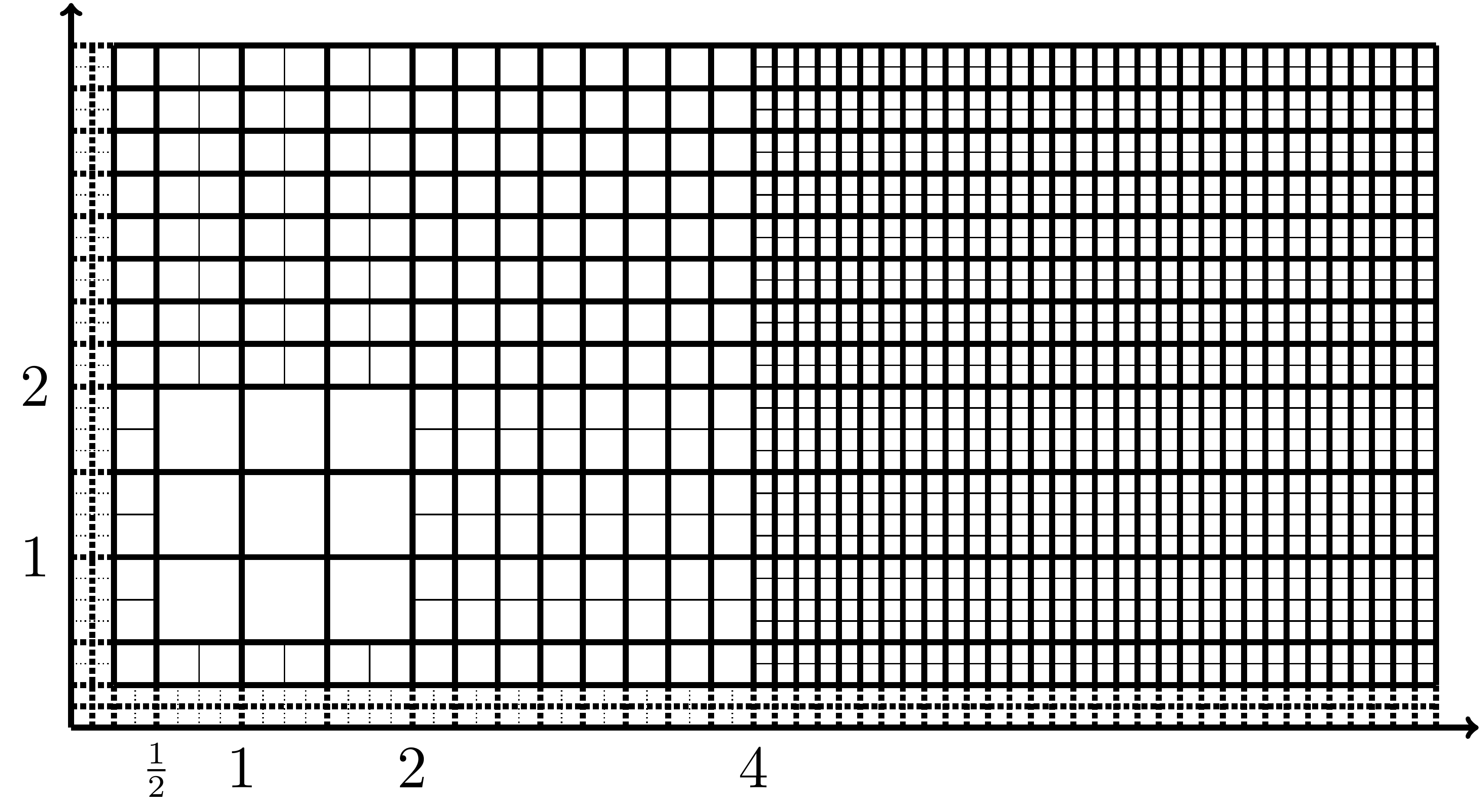}
\label{lag1_figure}
\end{figure}

The one-dimensional version of $H^1\left(\LLag\right)$ was studied in \cite{JD-ConAppr-2008}. The admissible covering is the following
$$\qq_{L}=\set{[2^n+k2^{-n-1},2^n+(k+1)2^{-n-1}]\colon k=0,\ldots,2^{2n+1}-1, n\in\NN }\cup\set{[2^{-n},2^{-n+1}]\colon n\in \NN_+},$$
see Figure \ref{lag1_figure}. Using methods similar to those in \cite{JD-ConAppr-2008} we verify \eqref{gauss}--\,\eqref{a22} in Subsection \ref{ssec43}.\\


\cor{coro_lag}{
Let $\a_1,...,\a_d>-1/2$ and $L_{L} = L_{L}^{[\a_1]}+...+ L_{L}^{[\a_d]} $, be the multidimensional Laguerre operator on $L^2((0,\8)^d)$. Then, the Hardy spaces $H^1(L_{L})$, $H^1(L_{L}^\nu)$, $\nu\in(0,1)$, and $H^1_{at}(\qq_{L} \boxtimes ... \boxtimes \qq_{L})$ coincide. Moreover, the associated norms are comparable.
}

\subsubsection{ {\bf Schr\"odinger operators.}}\label{ssec143}
Let $L_S = - \Delta + V$ denote a Schr\"odinger operator on $\Rd$, where $V \in L^1_{loc}(\Rd)$ is a nonnegative potential. Since $V\geq 0$, we have
\eq{ \label{schrodinger-a0}
0\leq T_{S,t}(x,y) \leq \Ht(x,y),\qquad x,y \in \RR^d, t>0,}
where $T_{S,t} = \exp(-tL_S)$ and $H_t=\exp(t\Delta)$, see \eqref{heat_kernel}.
Following \cite{DZ_Studia_DK}, for fixed $V$, we assume that there is an admissible covering $\qq_S$ of $\Rd$ that satisfies the following conditions:
there exist constants $\rho>1$ and $\sigma > 0$ such that
\eq{\label{schrodinger-d} \tag{D'}
\sup_{y\in \s{Q}{1}} \int_{\Rd} T_{S,2^n d_Q^2}(x,y) \, dx \leq C \rho^{-n}, \quad Q\in\qq_S, \, n\in \NN,
}
\eq{\label{schrodinger-k} \tag{K}
\sup_{y\in \Rd} \int_0^{t} \int_{\Rd} H_s(x,y)\chi_{\s{Q}{3}}(x) V(x) \, dx \, ds \leq C \left( \frac{t}{d_Q^2} \right)^\sigma, \quad Q\in\qq_S, \, t\leq d_Q^2.
}

The Hardy spaces related to Schr\"odinger operators have been widely studied. It appears that for some potentials the atoms for $H^1(L_S)$ have local nature (as in our paper), but this is no longer true for other potentials. The interested reader is referred to \cite{DZ_Studia_DK,Dziubanski2017,DZ_Revista2,DZ_Annali,DP_Potential,DZ_Potential_2014,Hofmann_Memoirs,CzajaZienkiewicz_ProcAMS}.

In \cite{DZ_Studia_DK} the authors study potentials as above, but instead of assuming \eqref{schrodinger-d} they have a bit more general assumption $(D)$, which instead of $\rho^{-n}$ has an arbitrary summable sequence $(1+n)^{-1-\e}$ on the right-hand side of \eqref{schrodinger-d}. Moreover, the assumptions \eqref{schrodinger-d} and \eqref{schrodinger-k} are easy to generalize for products, see \cite[Rem. 1.8]{DP_Potential}. Therefore, for Schr\"odinger operators Theorem \ref{mainthm} is a bit weaker than results of \cite{DZ_Studia_DK}. However, Theorem \ref{mainthm2} gives additionally characterization by the semigroups $\exp\eee{-t{L_S^\nu}}$, $0<\nu<1$, provided that the stronger assumption \eqref{schrodinger-d} is satisfied. Let us notice that indeed \eqref{schrodinger-d} is true for many examples, including $L_S$ in dimension one with any nonnegative $V\in L^1_{loc}(\RR)$, see \cite{CzajaZienkiewicz_ProcAMS}.

In Subsection \ref{ssec43} we prove that \eqref{schrodinger-d} and \eqref{schrodinger-k} imply the assumptions of Theorems \ref{mainthm} and \ref{mainthm2}, which leads to the following.

\cor{coro_Schro}{
Let $L_S$ be given with a nonnegative $V \in L^1_{loc}(\Rd)$ and an admissible covering $\qq_S$ of $\Rd$. Assume that \eqref{schrodinger-d} and \eqref{schrodinger-k} are satisfied. Then the spaces
$H^1(L_S)$, $H^1(L_S^\nu)$, $\nu\in(0,1)$, and $H^1_{at}(\qq_S)$ coincide and the corresponding norms are equivalent.
}

\subsubsection{ {\bf Product of local and nonlocal atomic Hardy space.}}\label{sssec13555}
As we have mentioned, all atoms on the Hardy space $H^1(\RR^{d_1})$  satisfy cancellation condition, i.e. they are nonlocal atoms. However, if we consider the product $\Rd = \RR^{d_1} \times \RR^{d_2}$ and the operator $L= -\Delta + L_2$, where $L_2$  and $\qq_2$ satisfies the assumptions \eqref{gauss}--\,\eqref{a22} on $\RR^{d_2}$ then the resulting Hardy space $H^1(L)$ shall have local character.

More precisely, if $\RR^{d_1}\boxtimes \qq_2$ is the admissible covering that arise by splitting all the strips $\RR^{d_1} \times Q_2$, $Q_2 \in \qq_2$, into countable many cuboids $Q_{1,n} \times Q_2$, where $Q_{1,n} = Q(z_n, d_{Q_2})$. Then we have the following corollary (see Subsection \ref{ssec44}).

\cor{coro_locnonloc}{
Let $L=-\Delta + L_2$, where $-\Delta$ is the standard Laplacian on $\RR^{d_1}$ and $L_2$ with an admissible covering $\qq_2$ of $\RR^{d_2}$ satisfy \eqref{gauss}--\,\eqref{a22}. Then the spaces
$H^1(L)$, $H^1(L^\nu)$, $\nu\in(0,1)$, and $H^1_{at}(\RR^{d_1}\boxtimes \qq_2)$ coincide and the corresponding norms are equivalent.
}

\subsection{Organization of the paper}\label{ssec14}

The paper is organized in the following way. Section \ref{sec2} is devoted to prove some preliminary estimates and to recall some known facts about local Hardy spaces on $\Rd$. In Section \ref{sec3} we prove our main results, namely Theorems \ref{atom-hardy}, \ref{mainthm}, and \ref{mainthm2}. In Section \ref{sec4} we prove that the examples given in Subsection \ref{ssec13} satisfy assumptions \eqref{gauss}--\,\eqref{a22}. We use standard notation, i.e. $C$ denotes some constant that can change from line to line.

\section{Preliminaries}\label{sec_local}\label{sec2}

\subsection{Auxiliary  estimates.}\label{ssec31}

For an admissible covering $\qq$ of {$X$ let us denote for $Q\in \qq$ the functions $\psi_Q\in C^1(X)$  satisfying
\eq{\label{partition}
0 \leq \psi_Q (x) \leq \chi_{Q^*}(x), \quad \norm{\psi_Q'}_\8 \leq C d_Q^{-1}, \quad \sum_{Q\in \qq} \psi_Q(x) = \chi_X(x).
}

It is easy to observe that such family $\set{\psi_Q}_{Q\in \qq}$ exists, provided that $\qq$ satisfies Definition \ref{finite_covering}.  The family $\set{\psi_Q}_{Q\in\qq}$ shall be called  {\it a partition of unity} related to $\qq$.

\prop{a3a4}{
Assume that $T_t$, and an admissible covering $\qq$ satisfy \eqref{a0} and \eqref{a1}. Let $\psi_Q$ be a partition of unity related to $\qq$. Then
\eq{ \label{a3}
\sup_{y \in \s{Q}{1}} \int_{\s{Q}{2}} \sup_{t > d_Q^2} \Tt(x,y) dx \leq C, \quad Q\in \qq,
}
and
\eq{ \label{a4}
\sup_{y\in X} \sum_{Q\in\qq} \int_{\s{Q}{2}}  \sup_{t\leq d_Q^2} \Tt(x,y) \abs{\psi_Q(x) - \psi_Q(y)} dx \leq C.
}
}

\pr{
By \eqref{a0} we have $\Tt(x,y) \leq C t^{-d/2}$. Obviously, $|Q^{**}| \leq C|Q| \leq C d_Q^{d}$, hence
\eqx{
\sup_{y\in\s{Q}{1}} \int_{\s{Q}{2}} \sup_{t>d_Q^2} \Tt(x,y) \, dx \leq  \int_{\s{Q}{2}} \sup_{t>d_Q^2} t^{-d/2} \, dx \leq C.}

We now turn to prove \eqref{a4}. Fix $y\in X$ and $Q_0\in \qq$  such that $y\in Q_0$ . Denote $N(Q_0)= \set{Q\in\qq \ : \ \s{Q_0}{3} \cap \s{Q}{3} \neq \emptyset }$ (the neighbors of $Q_0$) . Notice that $|N(Q_0)|\leq C$, see \eqref{finite_covering}. Then
\spx{
\sum_{Q \in \qq} \int_{\s{Q}{2}} \left[ \sup_{ t \leq d_Q^2  } \Tt(x,y) \abs{\psi_Q(x) - \psi_Q(y)} \right] \, dx = \sum_{Q \in N(Q_0)} ... + \sum_{Q \in \qq\setminus N(Q_0)}... \,  =: S_1 + S_2.}

Notice that for $Q \in N(Q_0)$ we have $d_Q \simeq d_{Q_0}$. To deal with $S_1$ we use \eqref{a0} and the mean value theorem for $\psi_Q$,
\spx{
  \sum_{Q \in N(Q_0)} \int_{\s{Q}{2}}  \sup_{  t \leq d_Q^2  } \Tt(x,y) \abs{\psi_Q(x) - \psi_Q(y)}  \, dx & \leq C \sum_{Q \in N(Q_0)} \int_{\s{Q}{2}} \ \sup_{ t>0 } \  t^\nu \left( t +|x-y|^2 \right)^{-d/2-\nu}
 \frac{|x-y|}{d_Q} \, dx\\
& \leq C \sum_{Q \in N(Q_0)} d_Q^{-1} \int_{\s{Q}{2}}  |x-y|^{-d+1} \, dx \\
& \leq C |N(Q_0)| d_{Q_0}^{-1} \int_{CQ_0} |x-y|^{-d+1}\, dx \leq C.
}

To estimate $S_2$ we use $\norm{\psi_Q}_{\8} \leq 1$ and \eqref{a1}, getting
\spx{
\sum_{Q \in \qq\setminus N(Q_0)} \int_{\s{Q}{2}}  \sup_{ t \leq d_Q^2 } \Tt(x,y) \abs{\psi_Q(x) - \psi_Q(y)}  \, dx & \leq 2 \sum_{Q \in \qq\setminus N(Q_0)} \int_{\s{Q}{2}}\sup_{t>0}  \Tt(x,y) \, dx\\
&  \leq C \int_{(\s{Q_0}{2})^c}  \sup_{t>0}  \Tt(x,y) \, dx \leq C.}
}

\lem{norm_ineq_lem}{Assume that $T_t$ satisfy \eqref{a0}. Then, for $f\in L^1(X) + L^\8(X)$,
\eqx{\norm{f}_{L^1(X)}\leq \norm{ \sup_{t>0}|T_t f|}_{L^1(X)}.}}

The proof of the Lemma \ref{norm_ineq_lem} goes by standard arguments. For the convenience of the reader we present details in Appendix.

\subsection{Local Hardy spaces.}\label{ssec22}
In this section, we recall some classical results on local Hardy spaces, see \cite{Goldberg_Duke}. Let $\tau>0$ be fixed. We are interested in decomposing into atoms a function $f$ such that
	\eq{\label{local_max}
    \norm{\sup_{t\leq \tau^2} \abs{H_t f}}_{L^1(\Rd)}<\8.
    }
It is known, that \eqref{local_max} holds if and only if $f(x) = \sum_k \la_k a_k(x)$, where $\sum_k \abs{\la_k} <\8$ and $a_k$ are either the classical atoms or the {\it local atoms at scale} $\tau$. The latter are atoms $a$ supported in a cube $Q$ of diameter at most $\tau$ such that $\norm{a}_\8 \leq |Q|^{-1}$ but we do not impose cancellation condition. In other words one may say that this is the space $H^1_{at}(\qq^{\{\tau\}})$ introduced in Subsection \ref{ssec12}, where $\qq^{\{\tau\}}$ is a covering of $\Rd$ by cubes with diameter $\tau$. The next proposition states the local atomic decomoposition theorem in a version that will be suitable for us in the proof of Theorem \ref{atom-hardy}. This proposition can be obtained by known methods from the global characterization of the classical Hardy space $H^1(\Rd)$. One may also check the assumptions from a general result of Uchiyama \cite[Cor. 1']{Uchiyama}. The details are left for the interested reader.

 \prop{prop_local}{
Let $\tau>0$ be fixed and $\wt{T}_t$ denote either $H_t$ or $P_{t^\nu, \nu}$, see \eqref{heat_kernel} and \eqref{poissone_kernel}. Then, there exists $C>0$ that does not depend on $\tau$ such that:
 \en{
 \item For every classical atom $a$ or atom of the form $a(x)=|Q|^{-1}\chi_{Q}(x)$, where $Q=Q(z,r_1,...,r_d)$ is such that $r_1\simeq ...\simeq r_d \simeq \tau$ we have
 	$$\norm{\sup_{t\leq \tau^2}\abs{\wt{T}_t a}}_{L^1(\Rd)} \leq C.$$

 \item If $f$ is such that $\supp f \subseteq Q^*$, where $Q=Q(z,r_1,...,r_d)$ is such that $r_1\simeq ...\simeq r_d \simeq \tau$, and
 	$$\norm{\sup_{t\leq \tau^2}\abs{\wt{T}_t f}}_{L^1(Q^*)}=M<\8,$$
then there exist sequences $\{\la_k\}_k$ and $\{a_k(x)\}_k$, such that $f(x) = \sum_k \la_k a_k(x)$, $\sum_k \abs{\la_k} \leq C M$, and $a_k$ are either the classical atoms supported in $Q^{*}$ or $a_k(x) = |Q|^{-1}\chi_{Q}(x)$.
 }
}

\rem{remark1}{
Proposition \ref{prop_local} remains 	valid for many other kernels $\wt{T}_t$ 	satisfying \eqref{a0} and, therefore, 	Theorem \ref{atom-hardy} holds for such kernels.
}


\section{Proofs of Theorems \ref{atom-hardy}, \ref{mainthm}, and \ref{mainthm2}.}\label{sec3}

\subsection{Proof of Theorem \ref{atom-hardy}}\label{ssec31}
\pr{
Recall that by the assumptions and Proposition \ref{a3a4} we also have  that \eqref{a3} and \eqref{a4} are satisfied. We shall prove two inclusions.

{\bf First inequality: $\norm{f}_{H^1(L)} \leq C \norm{f}_{H^1_{at}({\qq})}$.} It suffices to show that for every $\qq$-atom $a$ we have $\norm{\sup_{t>0} |\Tt a|}_{L^1(X)} \leq C$, where $C$ does not depend on $a$. Let $a$ be associated with a cuboid $Q\in\qq$, i.e. $\mathrm{supp}\,a \subset Q^*$. Recall that $\wt{T}_t$ is either $H_t$ or $P_{t^\nu, \nu}$, see \eqref{heat_kernel} and \eqref{poissone_kernel}. Observe that by using \eqref{a1}, \eqref{a2}, \eqref{a3}, and part {\bf 1.} of Proposition \ref{prop_local} we get

\spx{
	\norm{\sup_{t>0} |\Tt a|}_{L^1(X)}  \leq &\norm{\sup_{t>0} |\Tt a|}_{L^1((\s{Q}{2})^c)} + \norm{\sup_{t \leq d_Q^2} |(\Tt - \wt{T}_t )a|}_{L^1(\s{Q}{2})}\\
    & + \norm{\sup_{t> d_Q^2} |\Tt a|}_{L^1(\s{Q}{2})}  + \norm{\sup_{t \leq d_Q^2} |\wt{T}_t a|}_{L^1(\s{Q}{2})} \leq C.
}

{\bf Second inequality: $ \norm{f}_{H^1_{at}({\qq})} \leq C \norm{f}_{H^1(L)}$.} Assume that $\norm{\sup_{t>0} |\Tt f|}_{L^1(X)} <\infty$.  Let $\psi_Q$ be  a partition of unity related to $\qq$,  see \eqref{partition}. We have $f = \sum_{Q \in \qq} \psi_Q f$. Denote $f_Q = \psi_Q f$  and notice that since $\supp \, f_Q \subset \s{Q}{1}$, then
\eq{\label{ht-dec}
	\wt{T}_t  f_Q = (\wt{T}_t  - \Tt)f_Q + \left( {T}_t f_Q - \psi_Q\cdot\Tt f \right) + \psi_Q \cdot \Tt f.
}
Clearly,
\eq{ \label{psi-tt-fq}
	\sum_{Q\in\qq} \norm{ \sup_{t \leq d_Q^2} \abs{\psi_Q \Tt f}}_{L^1(\s{Q}{2})} \leq C \norm{\sup_{t>0} |\Tt f|}_{L^1(X)}.
}
Using \eqref{a2},
\eq{
	\sum_{Q\in\qq} \norm{ \sup_{t \leq d_Q^2} \abs{ (\wt{T}_t  - \Tt) f_Q}}_{L^1(\s{Q}{2})} \leq C \sum_{Q\in \qq} \norm{f_Q}_{L^1(X)} \leq C \norm{f}_{L^1(X)}.}
By \eqref{a4},
\sp{
\label{by-a4}
	\sum_{Q\in\qq} \norm{ \sup_{t \leq d_Q^2} \abs{\Tt f_Q - \psi_Q \cdot \Tt f }}_{L^1(\s{Q}{2})} &\leq \sum_{Q \in \qq} \int_X \abs{f(y)} \int_{\s{Q}{2}} \sup_{t\leq d_Q^2} \Tt(x,y) \abs{\psi_Q(y) - \psi_Q(x)} \, dx \, dy\\
    &\leq C \norm{f}_{L^1(X)}.
}
Using \eqref{ht-dec}--\eqref{by-a4} and Lemma \ref{norm_ineq_lem} we arrive at
\eqx{
	\sum_{Q\in\qq} \norm{ \sup_{t \leq d_Q^2} \abs{ \wt{T}_t  f_Q}}_{L^1(\s{Q}{2})} \leq C\norm{\sup_{t>0} |\Tt f|}_{L^1(X)}.
}
Now, from part {\bf 2.} of Proposition \ref{prop_local} for each $f_Q$ we obtain $\la_{Q,k}$, $a_{Q,k}$. Then
$$f = \sum_Q f_Q = \sum_{Q,k} \la_{Q,k} a_{Q,k}$$
and
$$\sum_Q \sum_k \abs{\la_{Q,k}} \leq C \sum_{Q\in\qq} \norm{ \sup_{t \leq d_Q^2} \abs{ \wt{T}_t  f_Q}}_{L^1(\s{Q}{2})} \leq C \norm{\sup_{t>0} T_t f}_{L^1(X)}. $$
Finally, we notice that all the atoms $a_{Q,k}$ obtained by Proposition \ref{prop_local} are indeed $\qq$-atoms.
}

\rem{remark111}{
The assumption \eqref{a0} has only been used in Proposition \ref{a3a4}. Therefore, in Theorem~\ref{atom-hardy} one may replace the assumption \eqref{a0} by the pair of assumptions \eqref{a3} and \eqref{a4}.
}

\subsection{Proof of Theorem \ref{mainthm}}\label{ssec32}
\pr{We shall show the following claim. If the assumptions \eqref{gauss}--\,\eqref{a22} hold for $\tti{j}_t(x_j,y_j)$ together with admissible coverings $\qq_j$ for $j=1,2$, then \eqref{gauss}--\,\eqref{a22} also hold for $\Tt(x,y) = \tti{1}_t(x_1,y_1)\cdot \tti{2}_t(x_2,y_2)$, together with $\qq = \qq_1 \boxtimes \qq_2$. This is enough, since by simple induction we shall get that in the general case $T_t(x,y) = \tti{1}_t(x_1,y_1)\cdot ... \cdot \tti{N}_t(x_N,y_N)$ with $\qq_1 \boxtimes ... \boxtimes \qq_N$ satisfy  \eqref{gauss}--\,\eqref{a22}, and, consequently, the assumptions of Theorem \ref{atom-hardy} will be fulfilled.

To prove the claim let $\tti{j}_t(x_j,y_j)$ and $\qq_j$ satisfy \eqref{gauss}--\,\eqref{a22} with $\gamma_j$ for $j=1,2$. Let $0<\gamma < \min(\gamma_1,\gamma_2)$ and fix $\delta \in [0,\gamma)$. Suppose that $\qq \ni Q \subseteq Q_1 \times Q_2$, where $Q_1 \in \qq_1$, $Q_2 \in \qq_2$, and without loss of generality we may assume that $d_{Q_1} \geq d_{Q_2}$. Hence, $Q=K\times Q_2$, where $K\subseteq Q_1$, see Definition \ref{def_box_prod} and Figure \ref{prostokat}. Denote by $z=(z_1,z_2)$ the center of $Q = K\times Q_2$. Obviously, \eqref{gauss} for the product follows from \eqref{gauss} for the factors.


\begin{figure}[h]
\caption{}\label{prostokat}
\vspace{1cm}
\centering
\includegraphics[scale=0.8]{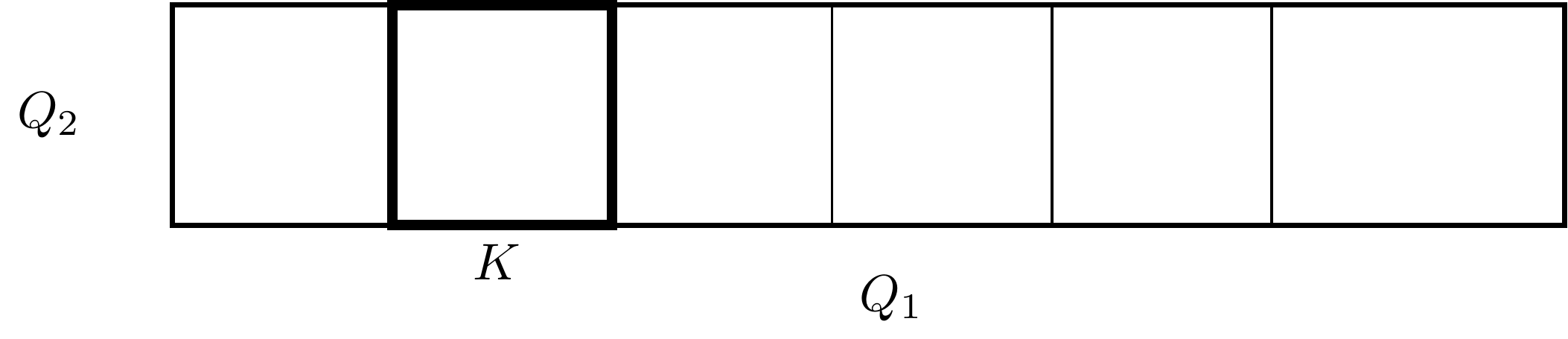}
\end{figure}

{\bf Proof of \eqref{a11} for $L_1+L_2$.}
Let  $y\in\s{Q}{1}$. Recall that $d_Q \simeq d_K \simeq d_{Q_2} \leq d_{Q_1}$. Let us write $(Q^{**})^c = S_1\cup S_2 \cup S_3$, where
$$S_1 = (K^{**})^c \times Q_2^{**}, \quad S_2 = K^{**} \times (Q_2^{**})^c, \quad S_3 =  (K^{**})^c \times (Q_2^{**})^c.$$

We start with $S_1$.
\spx{
\int_{S_1} \sup_{t>0} t^\delta \tti{1}_t(x_1,y_1) \tti{2}_t(x_2,y_2) \, dx  \leq &C \int_{(K^{**})^c} \sup_{t>0} {t^{-d_1/2-1/2}}\exp\left( -\frac{|x_1-y_1|^2}{ct} \right) \, dx_1\\
&\cdot \int_{Q_2^{**}} {\sup_{t>0} t^{-d_2/2+1/2+\delta}\exp\left( -\frac{|x_2-y_2|^2}{ct} \right)}\, dx_2 \\
\leq &  C \int_{(K^{**})^c} |x_1-z_1|^{-d_1-1} \, dx_1 \cdot \int_{Q_2^{**}} |x_2-z_2|^{-d_2+1+2\delta}\, dx_2  \\
\leq & C d_K^{-1} \cdot  d_{Q_2}^{1+2\delta} = C d_Q^{2\delta}.
}
The set $S_2$ is treated similarly. To estimate $S_3$ recall that $\delta < \gamma$. Using \eqref{gauss} for $\tti{1}_t(x_1,y_1)$ and \eqref{a11} for $\tti{2}_t(x_2,y_2)$ we arrive at
\spx{
	 \int_{S_3} \sup_{t>0} t^{\delta} \tti{1}_t(x_1,y_1) \tti{2}_t(x_2,y_2) \, dx \leq &C \int_{(K^{**})^c} \sup_{t>0} t^{-\gamma+\delta-d_1/2}\exp\left( -\frac{|x_1-y_1|^2}{ct} \right) \, dx_1\\&  \cdot \int_{(Q_2^{**})^c} \sup_{t>0} t^{\gamma} \tti{2}_t(x_2,y_2)\, dx_2 \\
	 \leq &C d_K^{-2\gamma+2\delta} d_{Q_2}^{2\gamma} \leq C d_Q^{2\delta}.}

{\bf Proof of \eqref{a22} for $L_1+L_2$.}
Let $y\in\s{Q}{1}$. In this proof $H_t$ is the classical heat semigroup on $\RR^{d_1}$, $\RR^{d_2}$ or on $\Rd$, depending on the context. First, notice that  by \eqref{gauss}, for constant $C>1$ and $i=1,2$, we have
\sp{\label{t-sim-dq}
\int_{\s{Q_i}{2}} \sup_{C^{-1} d_{Q_i}^2\leq t \leq Cd_{Q_i}^2} t^{-\gamma} \abs{T_t^{[i]}(x_i,y_i) - H_t(x_i,y_i)} \, dx_i
& \leq C d_{Q_i}^{-2\gamma} \int_{\s{Q_i}{2}}  d_{Q_i}^{-d_i} \exp\left( - \frac{|x_i-y_i|^2}{c d_{Q_i}^2} \right) \, dx_i\\ &\leq C d_{Q_i}^{-2\gamma}.
 }

Using the triangle inequality,
\spx{
& \int_{\s{Q}{2}} \sup_{t \leq d_Q^2 } t^{-\delta}\abs{\Tt(x,y) - \Ht(x,y)}dx \leq I_1 + I_2,
}
where
\alx{
I_1 &= \int_{\s{Q}{2}} \sup_{ t \leq d_Q^2 } t^{-\delta} \tti{1}_t(x_1,y_1) \abs{\tti{2}_t(x_2,y_2) - \Ht(x_2,y_2)} \, dx, \\
I_2 &=  \int_{\s{Q}{2}} \sup_{ t \leq d_Q^2  } t^{-\delta} \Ht(x_2,y_2) \abs{\tti{1}_t(x_1,y_1) - \Ht(x_1,y_1)}dx.
}

Applying \eqref{gauss} for $\tti{1}_t(x_1,y_1)$ and \eqref{a22} together with \eqref{t-sim-dq} for $\tti{2}_t(x_2,y_2)$,
\spx{
 I_1 & \leq C \int_{\s{K}{2}} \sup_{ t \leq d_Q^2  } t^{\gamma-\delta} \tti{1}_t(x_1,y_1) \, dx_1 \cdot \int_{Q_2^{**}} \sup_{ t \leq C d_{Q_2}^2 } t^{-\gamma} \abs{\tti{2}_t(x_2,y_2) - \Ht(x_2,y_2)} \, dx_2 \\
& \leq C d_K^{2\gamma-2\delta} d_{Q_2}^{-2\gamma} \simeq C d_Q^{-2\delta},
}
since $0 \leq \delta < \gamma < \min(\gamma_1,\gamma_2)$. Similarly, by \eqref{heat_kernel}, \eqref{a22}, and \eqref{t-sim-dq}, we have
\spx{
 I_2 & \leq C \int_{\s{Q_2}{2}} \sup_{ t \leq d_Q^2  } t^{\gamma-\delta} \Ht(x_2,y_2) \, dx_2 \cdot \int_{\s{Q_1}{2}} \sup_{ t \leq C d_{Q_1}^2 } t^{-\gamma} \abs{\tti{1}_t(x_1,y_1) - \Ht(x_1,y_1)} \, dx_1 \\
& \leq  C d_Q^{2\gamma-2\delta}d_{Q_1}^{-2\gamma} \leq C d_Q^{-2\delta},
}
since $d_{Q_1} \geq d_{Q_2} \simeq d_Q$.
}

\subsection{Proof of Theorem \ref{mainthm2}}\label{ssec33}

\pr{ For $\nu \in (0,1)$ the subordination formula introduced by Bochner \cite{SB_Fractional_PNASUSA} states that
\eq{\label{sub-pt}
P_{t^\nu,\nu}(x,y) = \int_0^\8 H_{ts}(x,y) d\mu_{\nu}(s),
}
and
\eq{\label{sub-kt}
K_{t^\nu,\nu}(x,y) = \int_0^\8 T_{ts}(x,y) d\mu_{\nu}(s),
}
where $\mu_\nu$ is a probability measure defined by the means of the Laplace transform $\exp(-x^\nu) = \int_0^\8 \exp(-xs) d\mu_\nu(s)$. By inverting the Laplace transform one obtains that $d\nu(s) = g_{\nu}(s)\,ds$ with
\eqx{\label{prob-density}
 0\leq g_\nu(s) = \int_0^\8 \exp\eee{ws\cos\theta_\nu+w^\nu\cos\theta_\nu}\sin\left(sw\sin\theta_\nu-w^\nu\sin\theta_\nu + \theta_\nu\right)  \, dw, \quad s>0,
}
where $\theta_\nu=\frac{\pi}{1+\nu}\in(\frac{\pi}{2},\pi)$, see \cite[Rem. 1]{KY_Fractional_PJA}. Notice that $\cos\theta_\nu <0$ and, therefore,
\sp{\label{density-est}
g_\nu(s) \leq \abs{\int_0^{s^{-1}}... \, dw}+\abs{\int_{s^{-1}}^\8... \, dw} \leq \int_0^{s^{-1}} \, dw + \int_{s^{-1}}^\8 \exp(ws\cos\theta_\nu) \, dw \leq C s^{-1}.
}
Assume that $T_t$ and $\qq$ satisfy \eqref{gauss}--\,\eqref{a22}. Then, Theorem \ref{mainthm2} follows from Theorem \ref{atom-hardy}, provided that we prove \eqref{a0}--\,\eqref{a2} for $K_{t^\nu,\nu}$ and $\qq$. First, notice that \eqref{a0} for $K_{t^\nu,\nu}$ follows from \eqref{sub-kt} and \eqref{gauss} for $T_t$. Coming to \eqref{a1}, let $Q\in\qq$ and $y\in\s{Q}{1}$. Since $\mu_\nu$ is a probability measure, using \eqref{sub-kt} and \eqref{a1} for $T_t$, we obtain
\spx{
	\int_{(Q^{**})^c} \sup_{t>0} K_{t^\nu,\nu}(x,y) \, dx & = \int_{(Q^{**})^c} \sup_{t>0} \int_0^\8 T_{st}(x,y) d\mu_\nu(s) \, dx\\
	& \leq  \int_0^\8 \int_{(Q^{**})^c} \sup_{t>0} T_{st}(x,y)  \, dx \, d\mu_\nu(s) \leq C.}
Having \eqref{a1} proved, we turn to \eqref{a2}. By \eqref{sub-pt}--\eqref{density-est}, and \eqref{a22} for $T_t$, we have
\spx{
\int_{Q^{**}} \sup_{t\leq d_Q^2} \abs{K_{t^\nu,\nu}(x,y) - P_{t^\nu,\nu}(x,y) } \, dx  = &\int_{Q^{**}} \sup_{t\leq d_Q^2} \abs{\int_{0}^\8 \left(T_{u}(x,y) - H_{u}(x,y)\right) \,  g_\nu(u/t) \, \frac{du}{t} } \, dx\\
 \leq& C \int_{Q^{**}} \sup_{t\leq d_Q^2} \int_{0}^\8 \abs{T_{u}(x,y) - H_{u}(x,y)}  (u/t)^{-1} \, \frac{du}{t}  \, dx \\
 \leq& C \int_{Q^{**}} \int_{0}^{d_Q^2} \abs{ T_{u}(x,y) - H_{u}(x,y)} \, \frac{du}{u}  \, dx \\
& + C \int_{Q^{**}} \int_{d_Q^2}^\8  \abs{ T_{u}(x,y) - H_{u}(x,y)} \, \frac{du}{u}  \, dx \\
	 \leq & C \int_{0}^{d_Q^2} u^{-1+\delta}  \int_{Q^{**}}  \sup_{u\leq d_Q^2} u^{-\delta}\abs{T_{u}(x,y) - H_{u}(x,y)}  \, dx  \, du \\
     &+ C \int_{Q^{**}} \int_{d_Q^2}^\8  u^{-d/2-1} \, du \, dx \\
     \leq &C  d_Q^{-2\delta} \int_{0}^{d_Q^2} u^{-1+\delta} \, du + Cd_Q^{d} d_Q^{-d} \leq C.
     }
This ends the proof of  Theorem \ref{mainthm2}.
}

\rem{rem_1234}{
It is worth to notice, that in the proof of \eqref{a2} for the subordinate semigroup $K_{t,\nu}$ we needed \eqref{a22} for $T_t$, not only \eqref{a2}.
}

\section{Applications}\label{sec4}
In this section for simplicity, we use the same notation $\Tt(x,y)$ for the integral kernels of semigroups generated by different operators.

\subsection{Bessel operator.}\label{ssec41}

Let us start with the following asymptotics of the Bessel function $I_\tau$,
\begin{align} \label{Bessel-function-small}
&&I_\tau(x) & = C_\tau x^\tau + O(x^{\tau+1}), & \text{ for } & x\sim 0,&&\\
\label{Bessel-function-large}
&&I_\tau(x) & = (2\pi x)^{-1/2} e^x + O(x^{-3/2}e^x), & \text{ for } & x \sim \8,&&
\end{align}
see e.g. \cite[p.\ 203--204]{Watson}.

\prop{a_bessel}{
Let $X=(0,\8)$ and $\be>0$. Then \eqref{gauss}--\,\eqref{a22} hold for $L_B^{[\be]}$ with $\qq_B$.
}
\pr{
We shall use similar ideas to those of \cite{BDT_d'Analyse}. The proof of \eqref{gauss} is well-known and follows almost directly from \eqref{bessel-kernel}, \eqref{Bessel-function-small} and \eqref{Bessel-function-large}. We skip the details. Let $\gamma \in(0,\min(1/2, \be/2))$ and $\delta\in[0,\gamma)$. Take $\qq_B \ni Q = [2^n,2^{n+1}]$, for some $n\in\ZZ$, and fix $y\in \s{Q}{1}$.

{\bf Proof of \eqref{a11}.} Notice that $y\simeq d_Q \simeq 2^n$. We have
\spx{
\int_{(\s{Q}{2})^c} \sup_{t>0} t^\delta \Tt(x,y) \, dx & \leq \int_0^\8 \sup_{t>xy} t^\delta \Tt(x,y) \, dx + \int_{(\s{Q}{2})^c} \sup_{t\leq xy} t^\delta \Tt(x,y) \, dx =: I_1 + I_2.}
Using \eqref{bessel-kernel} and \eqref{Bessel-function-small}, we obtain
\spx{
I_1 & \leq C \int_0^\8 \sup_{t>xy} (xy)^{\be} t^{\delta-\be-1/2} \exp\left( -\frac{x^2+y^2}{4t} \right)\, dx \\
& \leq C \int_0^\8 (xy)^{\be} (x^2+y^2)^{\delta-\be-1/2} \, dx \\
& = C y^{2\delta} \int_0^\8 x^\be (x^2+1)^{\de-\be-1/2}\, dx  \leq C d_Q^{2\de},
}
where in the last inequality we used the fact that $2\de<\be$.

Denote $z=3 \cdot 2^{n-1}$ (the center of $Q$). By  \eqref{bessel-kernel} and \eqref{Bessel-function-large},
\spx{
I_2 & \leq C \int_{(\s{Q}{2})^c} \sup_{t\leq xy} t^{\delta-1/2} \exp\left( -\frac{|x-y|^2}{4t} \right)\, dx \\
& \simeq  C \int_{(\s{Q}{2})^c} \sup_{t\leq xy} t^{\delta-1/2} \exp\left( -\frac{|x-z|^2}{{c}t} \right)\, dx \\
& \leq C \int_0^{2^n} \sup_{t>0} t^{\delta-1/2} \exp\left( -\frac{z^2}{{c_1}t} \right)\, dx + C \int_{2^{n+1}}^\8 \sup_{t\leq xy} t^{\delta-1/2} \exp\left( -\frac{x^2}{{c_2}t} \right)\, dx \\
&\leq C z^{2\de -1} 2^n + C \int_{2^{n+1}}^\8 (xy)^{\de-1/2} \exp\eee{-\frac{x}{{c_2}y}}\, dx\\
&\leq C d_Q^{2\de}.}

{\bf Proof of \eqref{a22}.} Now observe that if $y\in \s{Q}{1}$ and $x\in\s{Q}{2}$, then $x \simeq y \simeq d_Q$. Therefore, $\frac{xy}{2t} \geq c$, when $t\leq d_Q^2$. Using \eqref{bessel-kernel}, \eqref{Bessel-function-large}, and $\de < 1/2$, we arrive at
\spx{
\int_{\s{Q}{2}} \sup_{t\leq d_Q^2} t^{-\delta} \abs{\Tt(x,y) - \Ht(x,y)} \, dx &
\leq \int_{\s{Q}{2}} \frac{\sqrt{xy}}{2} \sup_{t\leq d_Q^2}  t^{-1-\de} \exp\left( -\frac{x^2+y^2}{4t} \right) \abs{I_{\be-\frac{1}{2}} \eee{\frac{xy}{2t}}-\frac{e^{\frac{xy}{2t}}}{\sqrt{\frac{\pi xy}{t}}}} \, dx\\
&\leq C \int_{\s{Q}{2}} \sup_{t\leq d_Q^2} t^{1/2-\delta} (xy)^{-1} \exp\left( -\frac{|x-y|^2}{4t} \right) \, dx \\
& \leq C d_Q^{1-2\delta} \cdot d_Q^{-2} \cdot d_Q \leq C d_Q^{-2\delta}.
}}
\subsection{Laguerre operator.}\label{ssec43}
Using the asymptotic estimates for the Bessel function \eqref{Bessel-function-small} and \eqref{Bessel-function-large} in formula \eqref{lag1-kernel}, one can obtain
\eq{\label{dfg}
T_t(x,y)\leq C t^{-1/2}\exp\left(-c\frac{|x-y|^2}{t}\right) e^{-ctxy} \min(1,(xy/t)^{\a+1/2}),\qquad x,y\in X, t>0,
}
see \cite[Eq. (2.12) and (2.13)]{JD-ConAppr-2008}.
\prop{a_lag}{
Let $X=(0,\8)$ and $\a>-1/2$. Then \eqref{gauss}--\,\eqref{a22} hold for $L_{L}^{[\a]}$ with $\qq_{L}$.
}

\pr{ We shall use similar estimates to those of \cite{JD-ConAppr-2008}. Note that \eqref{gauss} follows immediately from \eqref{dfg}. Let us fix positive constants $\gamma<\min(1/4,\a/2+1/4)$ and $\delta \in [0,\gamma)$. Fix $Q\in\qq_{L}$ and  $y\in\s{Q}{1}$.

{\bf Proof of \eqref{a11}.} We write
\eqx{
\int_{(\s{Q}{2})^c} \sup_{t>0} t^{\delta} T_t(x,y) \, dx =  \int_{(\s{Q}{2})^c\cap (0,d_Q)} ... +  \int_{(\s{Q}{2})^c\cap (d_Q, \8)} ... =: I_1 + I_2.
}
Since $|x-y|\geq C d_Q$ and $\delta < 1/2$, we have
\spx{
I_1 & \leq C \int_{(\s{Q}{2})^c\cap (0,d_Q)} \sup_{t>0} t^{\delta - 1/2} \exp\left( - \frac{|x-y|^2}{ct} \right) \, dx \\
& \leq C \int_{(\s{Q}{2})^c\cap (0,d_Q)} |x-y|^{2\delta - 1}  \, dx  \\
& \leq C d_Q^{2\delta-1} d_Q \leq C d_Q^{2\delta}.
}
In order to estimate $I_2$ we consider two cases depending on the localization of $Q$.\\

{\bf Case 1:} $Q=[2^{-n},2^{-n+1}]$, $n\in\NN_+$. In this case $y\simeq d_Q = 2^{-n}$. Observe that if $x\in (\s{Q}{2})^c\cap (d_Q, \8) $, then $|x-y|\sim x$ and
\spx{
\sup_{t>0}t^\delta T_t(x,y)& \leq C \sup_{t>0} t^{\delta -1/2}  \eee{\frac{xy}{t}}^{\a+1/2} \exp\left( - \frac{x^2}{ct} \right) \\
&\leq C d_Q^{\alpha+1/2}x^{2\delta-\alpha-3/2}.
}
Therefore,
$
I_2 \leq  C d_Q^{\alpha+1/2} \int_{d_Q}^{\8} x^{2\delta-\alpha-3/2} \, dx \leq C d_Q^{2\delta},
$
since $\delta \leq \a/2+1/4$.

{\bf Case 2:} $Q \subset [2^n,2^{n+1}]$, $n\in \NN$. Then $y^{-1} \simeq d_Q \simeq 2^{-n}$. Recall that $\de <1/2$. By using the inequality $\exp\eee{-cxyt} \leq C (xyt)^{-1}$ in \eqref{dfg}, we get
\spx{
I_2 & \leq C \int_{(\s{Q}{2})^c\cap (d_Q, \8)} \sup_{t>0} \ (xy)^{-1} t^{\delta -3/2} \exp\left( - \frac{|x-y|^2}{ct} \right) \, dx \\
& \leq C d_Q \int_{(\s{Q}{2})^c\cap (d_Q, \8)} x^{-1} |x-y|^{2\delta - 3} \, dx\\
& \leq C d_Q d_Q^{2\delta - 1 } \int_{(\s{Q}{2})^c\cap (d_Q, \8)} x^{-1} |x-y|^{-2} \, dx \\
& \leq Cd_Q^{2\delta} \left( \int_{(\s{Q}{2})^c\cap \left(d_Q,d_Q^{-1}/4\right)} d_Q^{-1} y^{-2} \, dx + \int_{(\s{Q}{2})^c\cap \left(d_Q^{-1}/4, \ \8\right)} d_Q |x-y|^{-2} \, dx  \right)  \leq C d_Q^{2\delta}.}

{\bf Proof of \eqref{a22}.}
For $x\in\s{Q}{2}$, $y\in\s{Q}{1}$ and $t\leq d_Q^2$, we apply an estimate that can be deduced from the proof of \cite[Prop.~2.3]{JD-ConAppr-2008}, namely
\eqx{
|T_t(x,y)-H_t(x,y)|\leq C t^{1/2}\eee{xy+(xy)^{-1}} \leq C t^{1/2} d_Q^{-2},
}
where the second inequality follows from the relation between $d_Q$ and the center of $Q$. Thus, for $\delta < 1/2$,
\eqx{
	\int_{Q^{**}}\sup_{t<d_Q^2} t^{-{\delta}} |T_t(x,y)-H_t(x,y)|\,dx\leq C d_Q^{-2}\int_{Q^{**}}\sup_{t<d_Q^2} t^{1/2-{\delta}}\,dx\leq C d_Q^{-2{\delta}}.
}
}

\subsection{Schr\"odinger operator.}\label{ssec42}

This subsection is devoted to proving the following proposition.
\prop{a_schrodinger}{
Let $L_S=-\Delta +V$ be a Schr\"odinger operator with $0\leq V \in L^1_{loc}(\Rd)$. Assume that for some admissible covering $\qq_S$ the conditions \eqref{schrodinger-d} and \eqref{schrodinger-k} hold. Then \eqref{gauss}--\,\eqref{a22} are satisfied for $L_S$ and $\qq_S$.
}

\pr{
In the proof we use estimates similar to those in \cite{DZ_Studia_DK}. For the completeness we present all the details. As we have already mentioned in \eqref{schrodinger-a0}, \eqref{gauss} holds since $V\geq 0$. Let us fix a positive $\gamma < \min(\log_2 \rho, \, \sigma)$, where $\rho$ and $\sigma$ are as in \eqref{schrodinger-d} and \eqref{schrodinger-k}, see subsection \ref{ssec143}. Consider $Q \in \qq_S$, $\delta\in[0,\gamma)$, and $y\in Q^{*}$.\\
{\bf Proof of \eqref{a11}.}  We have that
\spx{
\int_{(\s{Q}{2})^c} \sup_{t>0} t^\delta \, \Tt(x,y) \, dx & \leq \int_{(\s{Q}{2})^c} \sup_{t\leq 4d_Q^2} t^{\delta} \, \Tt(x,y) \, dx + \sum_{n\geq 2} \int_X \sup_{2^{n} d_Q^2 < t \leq 2^{n+1} d_Q^2} t^{\delta} \, \Tt(x,y) \, dx\\
&=:I_1 + I_2.
}

Denote by $z$ the center of cube $Q$. For $y\in Q^*$ and $x\not\in Q^{**}$ we have $d_Q \leq C |x-y| \simeq |x-z|$. Using \eqref{gauss} we obtain that
\spx{
I_1 &\leq C \int_{(\s{Q}{2})^c} \sup_{t\leq 4 d_Q^2} t^{-d/2+\delta} \exp\left( -\frac{|x-z|^2}{ct}\right) \, dx\\
&\leq C \int_{(\s{Q}{2})^c}  d_Q^{-d+2\delta} \exp\left( -\frac{|x-z|^2}{c \, d_Q^2}\right) \, dx \leq C d_Q^{2\delta}.
}
By \eqref{gauss} and \eqref{schrodinger-d},
\spx{
I_2 & \leq \sum_{n\geq 2} \int_{\Rd} \int_{\Rd} \sup_{2^{n} d_Q^2 < t \leq 2^{n+1} d_Q^2} t^{\delta} T_{t-2^{n-1} d_Q^2}(x,u) T_{2^{n-1} d_Q^2}(u,y) \, du \, dx\\
& \leq C \sum_{n\geq 1} (2^n d_Q^2)^{\delta}\int_{\Rd} T_{2^n d_Q^2}(u,y) \underbrace{\int_{\Rd}  (2^n d_Q^2)^{-d/2} \exp\left( -\frac{|x-u|^2}{c 2^n d_Q^2}\right)  \, dx}_{\leq C} \, du \\
& \leq  C d_Q^{2\delta} \sum_{n\geq 1} 2^{\delta n} \rho^{-n} \leq C  d_Q^{2\delta},
}
where in the last inequality we have used that $2^\delta  < \rho$.

{\bf Proof of \eqref{a22}.}
As in \cite[Lem. 3.11]{DZ_Studia_DK} we write $ V = \chi_{\s{Q}{3}} V + \chi_{(\s{Q}{3})^c} V =: V' + V''$. The perturbation formula states that
$H_t(x,y) - T_t(x,y) = \int_0^t \int_{\Rd} H_{t-s}(x,u) V(u) T_s (u,y)\, du\, ds,$ so
\spx{
t^{-\delta} \abs{\Ht(x,y) - \Tt(x,y)}   =  &t^{-\delta}  \int_{\Rd} \int_0^{t} H_{t-s}(x,u) V''(u) T_s(u,y) \, ds \, du\\
&+t^{-\delta}  \int_{\Rd} \int_0^{t/2} H_{t-s}(x,u) V'(u) T_s(u,y) \, ds \, du \\  & + t^{-\delta}  \int_{\Rd} \int_{t/2}^t H_{t-s}(x,u) V'(u) T_s(u,y) \, ds \, du \\
 =: &I_3(x,y) + I_4(x,y) + I_5(x,y).}
For $0< s < t \leq d_Q^2$, $x\in\s{Q}{2}$, $u\in(\s{Q}{3})^c$, we have that $d_Q \leq C |x-u|$ and
\eqx{t^{-\delta} H_{t-s}(x,u) \leq (t-s)^{-\delta} H_{t-s}(x,u) \leq C d_Q^{-d-2\delta}\exp\left( -\frac{|x-u|^2}{c \, d_Q^2} \right) }
and, consequently,
\spx{
\int_{\s{Q}{2}} \sup_{t\leq d_Q^2} I_3(x,y) \, dx & \leq C \int_{\s{Q}{2}} \int_{\Rd} \int_0^\8 d_Q^{-d-2\delta}\exp\left( -\frac{|x-u|^2}{c \ d_Q^2} \right) V''(u) T_s(u,y) \, ds \, du \, dx \\
& \leq C d_Q^{-2\delta} \int_{\Rd} \int_0^\8  V''(u) T_s(u,y) \, ds \, dz\\
& \leq C d_Q^{-2\delta}.
}
In the last inequality we have used equivalent form of \cite[Lem. 3.10]{DZ_Studia_DK}. To estimate $I_4$, denote $t_j = 2^{-j} d_Q^2$ for $j\geq 1$. Notice that
\sp{\label{I-j}
I_{4,j}(x,y) := \sup_{ t_j \leq t \leq t_{j-1}} I_4(x,y) & \leq C \sup_{ t_j \leq t \leq t_{j-1}} \int_{\Rd} \int_0^{t/2}  (t-s)^{-\delta} H_{t-s}(x,u) V'(u) T_s(u,y) \, ds \, du \\ & \leq C \int_0^{t_j} \int_{\Rd} t_j^{-d-\delta} \exp\left( -\frac{|x-u|^2}{c \, t_j} \right) V'(u) H_s(u,y) \, du \, ds.
}
Using \eqref{I-j} and then applying \eqref{schrodinger-k} we obtain
\spx{
\int_{\s{Q}{2}} \sup_{t\leq d_Q^2} I_4(x,y) \, dx & \leq \sum_{j\geq 1} \int_{\Rd} \sup_{t_j \leq t \leq t_j} I_{4,j}(x,y) \, dx \\
& \leq C \sum_{j\geq 1} t_j^{-\delta} \int_{\Rd} \int_0^{t_j} \underbrace{\int_{\Rd} t_j^{-d} \exp\left( -\frac{|x-u|^2}{c \, t_j} \right) \, dx}_{\leq C} V'(u) H_s(u,y)  \, ds \, du \\
& \leq C d_Q^{-2\delta} \sum_{j\geq 1} 2^{j\delta} \left( \frac{t_j}{d_Q^2} \right)^{\sigma} \leq C d_Q^{-2\delta} \sum_{j\geq 1} 2^{-j(\sigma-\delta)} \leq C d_Q^{-2\delta},}
since $\delta < \sigma$. Finally, $I_5(x,y)$ can be estimated by a similar argument. We skip the details.
}

\subsection{Products of local and nonlocal atomic Hardy spaces.}\label{ssec44}
In this section we consider operator $L=-\Delta+L_2$, where $-\Delta$ is the standard Laplacian on $\RR^{d_1}$ and $L_2$ together with an admissible covering $\qq_2$ of $X_2 \subseteq \RR^{d_2}$ satisfies \eqref{gauss}--\,\eqref{a22}. Obviously, the kernel of $\exp\eee{-tL}$ is given by $T_t(x,y) = H_t(x_1,y_1) \cdot T_t^{[2]}(x_2,y_2)$, where $x=(x_1,x_2) \in \RR^{d_1}\times X_2 \subseteq \RR^{d_1} \times \RR^{d_2} = \Rd$. One immediately see that  $\Tt(x,y)$  satisfies \eqref{gauss}. Moreover, almost identical argument to the proof of Theorem \ref{mainthm} shows that $\Tt$ with $\qq = \Rd \boxtimes \qq_2$ satisfies \eqref{a11} and \eqref{a22}. The details are left to the interested reader.


\section*{Appendix}
This appendix is devoted to prove Lemma \ref{norm_ineq_lem}. This proof uses standard methods, see e.g. \cite{Preisner_Studia}. We present details for the sake of completeness. In fact we prove a more general Proposition \ref{ae_conv}, from which Lemma \ref{norm_ineq_lem} follows immediately. Recall that we consider a semigroup of operators $T_t$ that is strongly continuous on $L^2(X)$ and has integral kernel $T_t(x,y)$ satisfying \eqref{a0}. We start with the following lemma.
\lem{appendix_lem_pom}{
    Suppose that $T_t$ satisfies \eqref{a0}. There exists a sequence $\set{t_n}_n$ such that $t_n \to 0$ and for every $r>0$ we have:
    \eq{\label{lim_to_0}
    \lim_{n\to \8} \int_{|x-y|>r} T_{t_n}(x,y) \, dy = 0,
    }
    \eq{\label{lim_to_1}
    	\lim_{n\to \8} \int_{|x-y|\leq r} T_{t_n}(x,y) \, dy = 1,
    }
for a.e. $x\in X$.
}
\pr{
	Let $\nu\in (0,1)$ be the constant from \eqref{a0}. Observe that
	\eqx{
	\int_{|x-y|>r} T_t(x,y) \, dy \leq C \int_{|x-y|>r} \frac{t^\nu}{(t+|x-y|^2)^{\frac{d}{2}+\nu}} \, dy = C \int_{|y| > \frac{r}{\sqrt{t}}} (1+|y|^2)^{-d/2-\nu} \, dy \rightarrow 0,
	}
	as $t\to 0$, and \eqref{lim_to_0} is proved (for every $\set{t_n}_n$ such that $t_n \to 0$).

To show \eqref{lim_to_1} observe that for $f\in L^2(X)$ we have $\lim_{t \to 0} T_{t} f$ converges to $f$ in $L^2(X)$, so we can choose a sequence with a.e. convergence. Applying this to functions $f_n(x) = \chi_{Q(0,n)}(x)$ and using a diagonal argument we obtain a sequence $\set{t_n}_n$, which goes to $0$, and such that for a.e. $x\in X$ we have
	\eq{\label{rtyu}
	\lim_{n\to \8} \int_{X} T_{t_{n}}(x,y) \, dy = 1.
    }
Thus, \eqref{lim_to_1} follows from \eqref{rtyu} and \eqref{lim_to_0}.
	}

\prop{ae_conv}{Assume that $T_t$ satisfies \eqref{a0} and let $f\in L^1(X)+L^\8(X)$. There exists a~sequence $\set{t_n}_n$ such that $t_n\to0$ and for almost every $x\in X$,
\eqx{
	\lim_{n\to \8} T_{t_n} f(x) = f(x).
}}
}

\pr{Let $\set{t_n}_n$ be the sequence from Lemma \ref{appendix_lem_pom}. By the Lebesgue differentiation theorem we have
\eq{\label{leb_diff_thm}
\lim_{s\to 0} \abs{Q(x,s)}^{-1} \int_{Q(x,s)} \abs{f(y) - f(x)} \, dy = 0}
	for almost every $x\in X$, since $f\in L^1(X)+L^\8(X) \subset L^1_{\rm{loc}}(X)$. Consider the set $A$ of points $x\in X$ such that we have \eqref{leb_diff_thm}, and, additionally, \eqref{lim_to_0}--\eqref{lim_to_1} hold for all rational $r > 0$. Obviously, such set has full measure. Fix $\e > 0$ and $x\in A$. We will show that $\abs{T_{t_n} f(x) - f(x)} \leq C\e$ for large $n\in\NN$. Let $r>0$ be a fixed rational number such that for $s<r$ we have
	\eq{\label{diff2}
	\int_{Q(x,s)} \abs{f(y) - f(x)} \, dy \leq \e \abs{Q(x,s)}.}
	Assume that $\sqrt{t_n} < r$ for large $n$. Write
	\spx{
	T_{t_n}f(x) - f(x) = & f(x) \left( \int_{\abs{x-y}\leq r} T_{t_n}(x,y) \, dy - 1\right) + \int_{|x-y|>r} T_{t_n}(x,y) f(y) \, dy \\
 & + \int_{|x-y| <\sqrt{t_n}} T_{t_n}(x,y) \left(f(y) - f(x)\right) \, dy \\
 & + \int_{\sqrt{t_n} \leq |x-y| \leq r} T_{t_n}(x,y) \left(f(y) - f(x)\right) \, dy\\
 & =: I_1 + I_2 + I_3 + I_4.
	}
	Applying \eqref{lim_to_1} we obtain that $\abs{I_1} < \e$ for $n$ large enough. To treat $I_2$ we consider two cases.\\
	{\bf Case 1: $f\in L^\8$.} Using \eqref{lim_to_0} we have that $\abs{I_2} <\e$ for $n$ large enough.\\
	{\bf Case 2: $f\in L^1$.} By \eqref{a0},
	\spx{
	\abs{I_2}  \leq C \int_{|x-y| > r} \frac{t_n^\nu}{(t_n+|x-y|^2)^{d/2+\nu}} \abs{f(y)} \, dy \leq C \frac{t_n^\nu}{(t_n+r^2)^{d/2+\nu}} \norm{f}_{L^1(X)} < \e,
	}
	for $t_n$ small enough. To estimate $I_3$ observe that $T_{t_n}(x,y) \leq C{t_n}^{-d/2}$ and $\abs{Q(x,\sqrt{t_n})} \simeq t_n^{d/2}$. Since $\sqrt{t_n} < r$, by applying \eqref{diff2} we obtain
	\eqx{
	\abs{I_3} \leq C t_n^{-d/2} \int_{|x-y| < \sqrt{t_n}} \abs{f(y) - f(x)} \, dy < C\e.}
	To deal with $I_4$ let $N=\lceil \log_2(r/\sqrt{t_n}) \rceil$, so that $r \leq \sqrt{t_n} 2^N \leq 2r.$ Define
	\eqx{S_k = \set{x\in X \ : \ r2^{-k} < |x-y| < r 2^{-k+1} }}
	 for $k=1,...,N$. Using \eqref{a0} and \eqref{diff2} we get
	\spx{
	\abs{I_4} & \leq C t_n^\nu \sum_{k=1}^N \int_{S_k} (t_n+|x-y|^2)^{-d/2-\nu} \abs{f(y) - f(x)} \, dy \\
	& \leq C t_n^{-d/2} \sum_{k=1}^N  (r2^{-k}/\sqrt{t_n})^{-d-2\nu} \int_{S_k} \abs{f(y) - f(x)} \, dy \\
	& \leq C \e t_n^{\nu} \sum_{k=1}^N  (r2^{-k})^{-d-2\nu} (r2^{-k})^{d} \\
	& \leq C \e (\sqrt{t_n} r^{-1}2^N)^{2\nu} \leq C \e.
	}
}


\medskip \medskip
\noindent
{\bf Acknowledgments:} The authors would like to thank Jacek Dziuba\'nski, B\l a\.zej Wr\'obel, and the reviewers for helpful comments.
\medskip

\bibliographystyle{amsplain}        

\begin{thebibliography}{10}

\bibitem{Auscher_unpublished}
P.~Auscher, X.T. Duong, and A.~McIntosh, \emph{Boundedness of banach space
  valued singular integral operators and hardy spaces}, Unpublished preprint,
  2005.

\bibitem{BDT_d'Analyse}
J.J. Betancor, J.~Dziuba\'nski, and J.L. Torrea, \emph{On {H}ardy spaces
  associated with {B}essel operators}, J. Anal. Math. \textbf{107} (2009),
  195--219.

\bibitem{SB_Fractional_PNASUSA}
S.~Bochner, \emph{Diffusion equation and stochastic processes}, Proc. Nat.
  Acad. Sci. U. S. A. \textbf{35} (1949), 368--370.

\bibitem{Coifman_Studia}
R.R. Coifman, \emph{A real variable characterization of {$H^{p}$}}, Studia
  Math. \textbf{51} (1974), 269--274.

\bibitem{CzajaZienkiewicz_ProcAMS}
W.~Czaja and J.~Zienkiewicz, \emph{Atomic characterization of the {H}ardy space
  {$H^1_L(\Bbb R)$} of one-dimensional {S}chr\"odinger operators with
  nonnegative potentials}, Proc. Amer. Math. Soc. \textbf{136} (2008), no.~1,
  89--94 (electronic).

\bibitem{Dziubanski_Houston}
J.~Dziuba\'nski, \emph{Hardy spaces associated with semigroups generated by
  {B}essel operators with potentials}, Houston J. Math. \textbf{34} (2008),
  no.~1, 205--234.

\bibitem{JD-ConAppr-2008}
\bysame, \emph{Hardy spaces for {L}aguerre expansions}, Constr. Approx.
  \textbf{27} (2008), no.~3, 269--287.

\bibitem{DP_Potential}
J.~Dziuba\'nski and M.~Preisner, \emph{On {R}iesz transforms characterization
  of {$H^1$} spaces associated with some {S}chr\"odinger operators}, Potential
  Anal. \textbf{35} (2011), no.~1, 39--50.

\bibitem{Dziubanski2017}
\bysame, \emph{Hardy spaces for semigroups with gaussian bounds}, Annali di
  Matematica Pura ed Applicata (1923 -) (2017).

\bibitem{DPRS}
J.~Dziuba\'nski, M.~Preisner, L.~Roncal, and P.~R. Stinga, \emph{Hardy spaces
  for {F}ourier-{B}essel expansions}, J. Anal. Math. \textbf{128} (2016),
  261--287.

\bibitem{DZ_Studia_DK}
J.~Dziuba\'nski and J.~Zienkiewicz, \emph{Hardy spaces {$H^1$} for
  {S}chr\"odinger operators with certain potentials}, Studia Math. \textbf{164}
  (2004), no.~1, 39--53.

\bibitem{DZ_Annali}
\bysame, \emph{Hardy spaces {$H^1$} for {S}chr\"odinger operators with
  compactly supported potentials}, Ann. Mat. Pura Appl. (4) \textbf{184}
  (2005), no.~3, 315--326.

\bibitem{DZ_Revista2}
\bysame, \emph{On {H}ardy spaces associated with certain {S}chr\"odinger
  operators in dimension 2}, Rev. Mat. Iberoam. \textbf{28} (2012), no.~4,
  1035--1060.

\bibitem{DZ_Potential_2014}
\bysame, \emph{A characterization of {H}ardy spaces associated with certain
  {S}chr\"odinger operators}, Potential Anal. \textbf{41} (2014), no.~3,
  917--930.

\bibitem{PG_WH_Studia}
P.~G{\l}owacki and W.~Hebisch, \emph{Pointwise estimates for densities of
  stable semigroups of measures}, Studia Math. \textbf{104} (1993), no.~3,
  243--258.

\bibitem{Goldberg_Duke}
D.~Goldberg, \emph{A local version of real {H}ardy spaces}, Duke Math. J.
  \textbf{46} (1979), no.~1, 27--42.

\bibitem{Hofmann_Memoirs}
S.~Hofmann, G.~Lu, D.~Mitrea, M.~Mitrea, and L.~Yan, \emph{Hardy spaces
  associated to non-negative self-adjoint operators satisfying
  {D}avies-{G}affney estimates}, Mem. Amer. Math. Soc. \textbf{214} (2011),
  no.~1007, vi+78.

\bibitem{MK_FCAA}
M.~Kwa\'{s}nicki, \emph{Ten equivalent definitions of the fractional {L}aplace
  operator}, Fract. Calc. Appl. Anal. \textbf{20} (2017), no.~1, 7--51.

\bibitem{Latter_Studia}
R.H. Latter, \emph{A characterization of {$H^{p}({\bf R}^{n})$} in terms of
  atoms}, Studia Math. \textbf{62} (1978), no.~1, 93--101.

\bibitem{Preisner_Studia}
M. Preisner, \emph{Atomic decompositions for {H}ardy spaces related to
  {S}chr\"{o}dinger operators}, Studia Math. \textbf{239} (2017), no.~2,
  101--122.

\bibitem{Song_Yan}
L.~Song and L.~Yan, \emph{A maximal function characterization for {H}ardy
  spaces associated to nonnegative self-adjoint operators satisfying {G}aussian
  estimates}, Adv. Math. \textbf{287} (2016), 463--484.

\bibitem{Stein}
E.M. Stein, \emph{Harmonic analysis: real-variable methods, orthogonality, and
  oscillatory integrals}, Princeton Mathematical Series, vol.~43, Princeton
  University Press, Princeton, NJ, 1993, With the assistance of Timothy S.
  Murphy, Monographs in Harmonic Analysis, III.

\bibitem{Uchiyama}
A.~Uchiyama, \emph{A maximal function characterization of {$H^{p}$} on the
  space of homogeneous type}, Trans. Amer. Math. Soc. \textbf{262} (1980),
  no.~2, 579--592.

\bibitem{Watson}
G.N. Watson, \emph{A treatise on the theory of {B}essel functions}, Cambridge
  Mathematical Library, Cambridge University Press, Cambridge, 1995, Reprint of
  the second (1944) edition.

\bibitem{KY_Fractional_PJA}
K.~Yosida, \emph{Fractional powers of infinitesimal generators and the
  analyticity of the semi-groups generated by them}, Proc. Japan Acad.
  \textbf{36} (1960), 86--89.

\end{thebibliography}

\def\cprime{$'$}
\providecommand{\bysame}{\leavevmode\hbox to3em{\hrulefill}\thinspace}
\providecommand{\MR}{\relax\ifhmode\unskip\space\fi MR }
\providecommand{\MRhref}[2]{%
  \href{http://www.ams.org/mathscinet-getitem?mr=#1}{#2}
}
\providecommand{\href}[2]{#2}

\end{document}